\documentclass{article}
\usepackage[a4paper, margin=1in]{geometry} 
\usepackage{natbib}
\usepackage{amsmath, amssymb} 
\usepackage{wrapfig}
\usepackage{amsmath}
\usepackage{hyperref}
\usepackage{booktabs}
\usepackage{rotating}
\usepackage{multirow}
\usepackage{array}
\usepackage{graphicx}
\usepackage{caption}
\usepackage{siunitx} 
\usepackage{titlesec}
\usepackage{subcaption} 
\usepackage{empheq}
\title{Enhancing PINN Performance Through Lie Symmetry Group}
\author{
    Ali Haider Shah\thanks{E-mail: alihaidershahhaider@gmail.com } \and
    Naveed R. Butt\thanks{E-mail: naveed.butt@giki.edu.pk
 } \and
    Asif Ahmad\thanks{E-mail: asifahmad7007@gmail.com (Corresponding author) } \and
    Muhammad Omer Bin Saeed\thanks{E-mail: omer.saeed@giki.edu.pk}
}
\date{}

\begin{document}

\maketitle

\begin{center}
\textsuperscript{Faculty of Engineering Sciences, Ghulam Ishaq Khan Institute of Engineering Sciences and Technology, Topi, Pakistan.}
\end{center}


\begin{abstract}
This paper presents intersection of Physics informed neural networks (PINNs) and  Lie symmetry group to enhance the accuracy and efficiency of solving partial differential equation (PDEs). Various methods have been developed to solve these equations. A Lie group is an efficient method that can lead to exact solutions for the PDEs that possessing Lie Symmetry. Leveraging the concept of infinitesimal generators from Lie symmetry group in a novel manner within PINN leads to significant improvements in solution of PDEs. In this study three distinct cases are discussed, each showing progressive improvements achieved through Lie symmetry modifications and adaptive techniques. State-of-the-art numerical methods are adopted for comparing the progressive PINN models. Numerical experiments demonstrate the key role of Lie symmetry in enhancing PINNs performance, emphasizing the importance of integrating abstract mathematical concepts into deep learning for addressing complex scientific problems adequately.
\end{abstract}


\section{Introduction}
Partial differential equations (PDEs) provide a framework for understanding intricate phenomena across multiple disciplines like physics, chemistry, and biology. As a result, discovering solutions to PDEs serves as a valuable and efficient method for investigating the dynamic behaviors inherent in these equations. Burgers' equation is a fundamental partial differential equation and convection–diffusion equation occurring in various areas of applied mathematics, such as fluid mechanics, gas dynamics, and traffic flow. It is one of famous quasi-linear partial differential equation. This kind of equation is related to Navier–Stokes equation with the pressure term removed \cite{liu2019lie}.
The Burgers equation was initially presented in a publication by Bateman, where he also presented a particular solution of the Equation. Subsequently, a notable sequence of papers spanning from 1939 to 1965 further explored Burger's equation \cite{ozics2003finite}. 
\\
Since it is conventional to use approaches linked to the intrinsic properties of PDEs for constructing analytic (exact) solutions, the classical Lie symmetry theory of PDEs emerges as one of the most inspiring source in this consideration \cite{zhang2023enforcing}. Lie groups are highly helpful for building conservation laws and determining the integrability, translational, and rotational invariance of PDEs \cite{bluman2008symmetry} \cite{olver1993applications}. A continuous Lie symmetry transformation in the context of PDEs refers to a continuous transformation connecting two solutions in the same set of PDEs. Naively, it maps the solution set of the same system to another solution set of that system. basic examples of Lie groups are translations and rotations etc.
Moreover, the invariant surface conditions produced by Lie symmetries are instrumental in determining invariant solutions of PDEs. These conditions serve as critical constraints that ensure the solutions remain unchanged under specific symmetry transformations. By imposing these conditions, it becomes possible to identify solutions that possess the desired symmetries, providing valuable insights into the behavior and properties of the underlying physical or mathematical systems described by the PDEs.
 PDEs then can be simplified by using symmetry reduction techniques, leading to equations with fewer independent variables. However, even with these approaches, the search for explicit precise solutions for some symmetry-reduced PDEs remains difficult \cite{turgut2017similarity}. Therefore Sophisticated numerical techniques are used to find the numerical solution of PDEs, such as the Least Square Cubic B-spline Finite Element Method (LS-CB-FEM)\cite{kutluay2004numerical}, Modified Cubic B-spline Differential Quadrature Method (MCB-DQM)\cite{arora2013numerical},  Weighted Average Differential Quadrature Method (WA-DQM)\cite{jiwari2013numerical}, Finite Element Method (FEM)\cite{ozics2003finite}, Modified cubic B-spline collocation method (MCB-CM) \cite{mittal2012numerical} and Reproducing kernel function \cite{xie2008numerical}. Specifically Burger's equation of different orders and degrees has been successfully addressed by using the aforementioned numerical methods.\\
 In recent years, Physics Informed Neural Networks (PINN) has earned significant attention through the groundbreaking work of Raissi. In his work, he utilized neural networks to find the solution of PDE's. This novel approach, which falls under the deep learning domain, has attracted a lot of interest due to its capacity to facilitate data-driven PDE solution finding \cite{raissi2017physics}. So far, PINN has found extensive application in scientific computing, specifically in addressing both forward and inverse problem in non linear PDEs \cite{pu2021solving,zhang2023enforcing,blechschmidt2021three,wang2021data,kumar2024approximate,blechschmidt2021three}.
\\
The approximate solution by PINN might not always be feasible and can show significant departures from the desired result. Observation suggests that a variety of reasons might be behind these differences. Neural networks incorporate many factors, and to obtain the best results, each parameter must be well-optimized. For this reason, we continuously modify and enhance our neural networks to get optimized solution. This path of development has led to the modified versions of PINN. For instance, to address both forward and inverse PDE problems, for example, a gradient-enhanced PINN has been developed, which incorporates gradient-based information from the PDEs residual into the loss function \cite{yu2022gradient}. Lin and Chen proposed
a two-stage PINN where the first stage is a conventional  standard PINN and the second stage is to combine the
conserved quantities with mean squared error loss to train neural networks \cite{lin2022two}.Zhu et al. built group equivariant neural networks that effectively simulated various periodic solutions of nonlinear dynamical lattices while adhering to spatiotemporal parity symmetries \cite{zhu2022neural}. To speed up training and convergence, Jagtag et al. presented local adaptive activation functions by adding a scaling parameter to each layer and neuron independently \cite{jagtap2020locally}. Some of useful techniques such as properly designed non-uniform training weighing techniques are helpful for improving the accuracy of PINN \cite{wang2021understanding}. To remove obstacles in more complicated and realistic applications, more focused neural networks are developed, including the Bayesian PINN \cite{yang2021b}, the discrete PINN structure based on graph convolutional network and variational structure of PDE \cite{gao2022physics}, extended physics-informed neural networks (XPINNs) \cite{jagtap2020extended}, peripheral physics-informed neural network (PPINN) \cite{meng2020ppinn}, auxiliary physics-informed neural networks (A-PINN) \cite{yuan2022pinn}, fractional physics-informed neural networks (fPINNs) \cite{pang2019fpinns}, and variational PINN \cite{kharazmi2021hp}. To solve PDEs, neural networks require independent information. Leveraging some of the intrinsic properties of PDEs through Lie symmetry in PINN can be advantageous for enhancing the solution. Embedding symmetry loss information in total neural network loss has produced good results. \cite{li2023utilizing,zhang2023enforcing}.
\\
\\
The main idea of this study is to enhance PINN performance by using LIE group symmetry and compare the PINN model with state-of-the-art numerical methods. Burger's equation is adopted for the numerical experiment. Burger's equation is a fundamental PDE that serves as a case study in various scientific disciplines. This research is divided into three cases. Four up-to-date activation functions are used in each case with the same structural components of PINN. 
In Case A, the conventional PINN model is used whereas, in the next case (case B), we introduce Lie group symmetry infinitesimal transformations in a novel manner. The loss function is modified by using these transformations leading to progressive improvements from case A. To address  discrepancies, The PINN model  after utilizing the LIE group in case B is pronounced as modified symmetry-based PINN (m-SPINN) 
Case C involves the execution of the adaptive activation function technique to further enhance the results of previous case B. After introducing  adaptivity in m-SPINN we refer to it as the modified adaptive Symmetry-based PINN (m-ASPINN)
\\
While all preceding cases are optimized using the Adam optimizer.
By Numerical experiment in section 3, it is clear that using Lie symmetry for PDEs plays an important role in enhancing the results and efficiency of PINN. These results emphasize how crucial it is to incorporate sophisticated mathematical ideas into machine learning algorithms to achieve significant improvements in accuracy and performance.

\section{Methodology}

We take the following one dimensional non-linear Burger's equation:
\begin{equation}
 \frac{\partial u}{\partial t} + a u \frac{\partial u}{\partial x} - v \frac{\partial^2 u}{\partial x^2} = 0, \quad \quad (x,t) \in X \times [0, T]  
 \label{eq1}
\end{equation}

\begin{equation}
 u(x,0)=f(x), \quad \quad   x \in X \subset (a,b)
  \label{eq2}
\end{equation}

\begin{equation}
u(a, t) = 0 \quad \text{and} \quad u(b, t) = 0, \quad t \in [0, T].
 \label{eq3}
\end{equation}

 These equations show the behavior of the solution \( u \) with respect to time \( t \) and spatial coordinate \( x \), where \( a \), \( b \), \( m \), and \( T \) are constants, and \( f(x) \) represents the initial value of \( u \).
Where \(v\) is a parameter having a small value that represents the coefficient of kinematic viscosity.
Since eqn \eqref{eq1} has been addressed by lots of numerical methods e.g., \cite{arora2013numerical} \cite{kutluay2004numerical}\cite{jiwari2013numerical} \cite{ozics2003finite} \cite{mittal2012numerical}
\cite{xie2008numerical}\cite{korkmaz2011quartic}\cite{chung2010asymptotic}\cite{salas2010symbolic} \cite{lin2010new}\cite{korkmaz2006numerical}\cite{korkmaz2010numerical}\cite{dag2005b}
\cite{xie2008numerical} and \cite{mittal2012numerical}\ Therefore, we find it simple to compare with our proposed technique. \\
\\
We examine a feed-forward Neural Network shown in Figure \ref{fig:PINN} to comprehend the behavior of a neural network. Consider a NN of \( L \) layers. This structure consists of a single input layer, \( L - 1 \) hidden layers and one output layer. \(N_{l}\) are the number of neurons for each \( l \)-th layer. In this network, each layer experiences an output \( Z_{l-1} \in \mathbb{R}^{N_{l-1}} \) from the previous layer, which goes through an affine transformation, followed by the application of non-linear activation function \( \sigma(\cdot) \).  defined mathematically in eqn \eqref{eq4}.

\begin{equation}
Z_{l} = \sigma(W_{l} Z_{l-1} + b_{l})
 \label{eq4}
\end{equation}

Here, \( W_{l} \in \mathbb{R}^{N_{l} \times N_{l-1}} \) are weights and \( b_{l} \in \mathbb{R}^{N_{l}} \) biases of the \( l \)-th layer, respectively. For the output layer, we use different activation functions according to our required result. Thus, the final output of the NN is given by eqn \eqref{eq5}:

\begin{equation}
 O_{\theta}(Z) = Z_{L}
  \label{eq5}
\end{equation}
 
where \( Z \) = \( Z_{0}\)  represents the input, \( \theta = \{ W_l, b_l \}_{l=1}^{L} \) are the trainable parameters of the NN and \( O \) signifies the output of the NN. Furthermore, we use different optimizers to find the best possible result according to the nature of our problem and defined loss function.
In general loss function can be defined as:

\begin{equation}
\text{loss} = \sum_{i=1}^{n} (y_i - O_{\theta}(Z_i))^2   
 \label{eq6}
\end{equation}

where \(y_i\) represents the actual value and \(O_{\theta}(Z_i)\) represents the predicted value by NN for the \(i\)-th sample, and \(n\) is the total number of sample points.This is a generalized structure to represent neural network its architecture changes according to the nature of the problem. In the next section 2.1, we explain the variant of the neural network called PINN that is used for the solution of PDEs.
\\
\subsection{PINN based-solution for Burgers equation}
\begin{figure}[htbp] 
    \centering
    \includegraphics[width=0.9\textwidth]{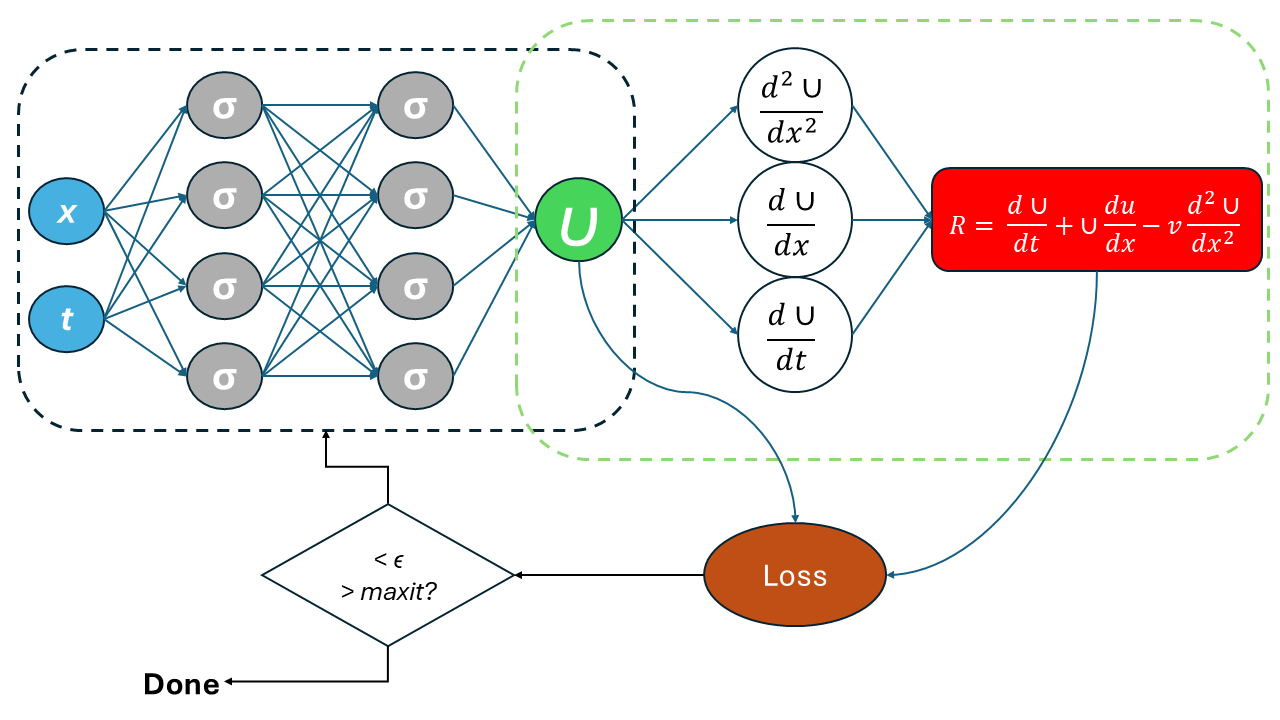} 
    \caption{Physics Informed Neural Network for Burgers' Equation}
    \label{fig:PINN} 
\end{figure}

The basic idea behind PINN for solving PDEs is to utilize a neural network $u_{\text{NN}}(t, x, \Theta)$ parameterized by $\Theta$ to approximate the exact solutions $u(t, x)$. In this technique, the neural network architecture is designed to encode the fundamental physics of the system, enabling it to learn the solution behavior directly from the data as shown in Figure\eqref{fig:PINN}.
\\
\\
Let's denote the data for the initial condition and boundary condition as \( \{(t_i^0, x_i^0, u_i^0)\}_{i=1}^{N_0} \) and \( \{(t_i^b, x_i^b, u_i^b)\}_{i=1}^{N_b} \), respectively, where \( N_0 \) and \( N_b \) represent the respective numbers of data points. Additionally, let \( \{(t_j^c, x_j^c)\}_{j=1}^{N_c} \) represent the collocation points used for evaluating the residual of the PDE, where \( N_c \) is the number of collocation points.
The loss term for the initial condition is denoted as \( \mathcal{L}_{\text{init}}(\Theta) \), for the boundary condition as \( \mathcal{L}_{\text{bound}}(\Theta) \), and for the residual term as \( \mathcal{L}_{\text{res}}(\Theta) \). These terms are defined as:

\begin{equation}
\mathcal{L}_{\text{init}}(\Theta) = \frac{1}{N_0} \sum_{i=1}^{N_0} \left( u_{\text{NN}}(t_i^0, x_i^0, \Theta) - u_i^0 \right)^2
 \label{eq7}
\end{equation}

\begin{equation}
\mathcal{L}_{\text{bound}}(\Theta) = \frac{1}{N_b} \sum_{i=1}^{N_b} \left( u_{\text{NN}}(t_i^b, x_i^b, \Theta) - u_i^b \right)^2
 \label{eq8}
\end{equation}

\begin{equation}  \mathcal{L}_{\text{res}}(\Theta) = \frac{1}{N_c} \sum_{j=1}^{N_c} \left( \mathcal{R}(t_j^c, x_j^c, \Theta) \right)^2  
\label{eq9}
\end{equation}

where $u_{\text{NN}}(t_i^0, x_i^0, \Theta)$ ,$u_{\text{NN}}(t_i^b, x_i^b, \Theta)$ represents the predicted solution of PDE at initial and boundary points by NN. The residual of the PDE at collocation points denoted as
$\mathcal{R}(t_j^c, x_j^c, \Theta)$.
Residual term for eqn \eqref{eq1} defined as:

\begin{equation}
\mathcal{R}= \frac{\partial u}{\partial t} + a u \frac{\partial u}{\partial x} - v \frac{\partial^2 u}{\partial x^2}
 \label{eq10}
\end{equation}

with differential operator $\mathfrak{R}$,

\begin{equation}
     \mathfrak{R} = \frac{\partial }{\partial t} + a u \frac{\partial }{\partial x} - v \frac{\partial^2 }{\partial x^2}   \label{eq11} 
\end{equation}

Finally  total loss function combining all terms in eqn \eqref{eq7}, eqn \eqref{eq8} and eqn \eqref{eq9} can be written as:

\begin{equation}
 \mathcal{L}(\Theta) = \alpha \mathcal{L}_{\text{init}}(\Theta) + \beta \mathcal{L}_{\text{bound}}(\Theta) + \gamma \mathcal{L}_{\text{res}}(\Theta)  
  \label{eq12}
\end{equation}

where $\alpha$, $\beta$ and $\gamma$ are weights, usually they are considered as $\alpha$ = $\beta$ = $\gamma$ = 1.
The PINN technique in this situation often makes use of the Adam \cite{kingma2014adam} or L-BFGS \cite{liu1989limited} algorithms to update the parameters $\Theta$ and, as a result, minimize the mean square error (MSE) loss function. We can put constraints that contain Independent properties and information of  our  system to make neural network more efficient. Thus, we extend this idea and include Lie symmetry group- based information of Burger's equation \eqref{eq1} into the loss function. Consequently, the PINN model converts into m-SPINN. So, in the next section, we discuss integration of PINN with Lie symmetry group for Burger's equation. 

\subsection{Solving the Burgers equation using m-SPINN}
One of PDE's intrinsic but hidden properties is symmetry. It offers a broadly applicable method for locating closed form solutions to PDEs, as the majority of practical systematic approaches, such as separation of variables, require the direct use of the symmetry method. System \eqref{eq1} admits Lie symmetry by a classical way of locating a local one-parameter  ($\epsilon$) Lie group with infinitesimal transformation. 
\\
The transformed variables are given by:

\begin{equation}
 \tilde{x} = x + \varepsilon \xi_1(x,t,u) + O(\varepsilon^2) 
  \label{eq13}
 \end{equation}   
 
 \begin{equation}
     \tilde{t} = t + \varepsilon \xi_2(x,t,u) + O(\varepsilon^2)    \label{eq14} 
 \end{equation}
 
\begin{equation}
       \tilde{u} = u + \varepsilon \eta(x,t,u) + O(\varepsilon^2)
        \label{eq15}
\end{equation}

This transformation leaves the system invariant.Lie’s third fundamental theorem tells that such a group is completely characterized by the infinitesimal generator.

\begin{equation}
\mathbf{L} = \xi_1\frac{\partial}{\partial x} + \xi_2\frac{\partial}{\partial t} + \eta\frac{\partial}{\partial u}  
 \label{eq16}
\end{equation}

we usually  pronounce the Lie group transformations described in equations \eqref{eq13}, \eqref{eq14}, \eqref{eq15} and the related operator $\mathbf{L}$ in eqn \eqref{eq16} as Lie symmetry without differentiating them.
The invariance condition are used to derive the infinitesimals $\xi_1(x, t, u)$, $\xi_2(x, t, u)$, and $\eta(x, t, u)$. For this we use second prolongation formula denoted as $\mathbf{L^{(2)}}$. we will neglect some non-useful term of $\mathbf{L^{(2)}}$  \cite{turgut2017similarity}. Second prolongation $\mathbf{L^{(2)}}$ for system \eqref{eq1} is defined as: 

\begin{equation}
 \mathbf{L^{(2)}} = \xi_1 \frac{\partial}{\partial x} + \xi_2 \frac{\partial}{\partial t} + \eta \frac{\partial}{\partial u} 
     + \pi_x \frac{\partial u}{\partial x} + \pi_t \frac{\partial u}{\partial t} + \pi_{xx} \frac{\partial^2 u}{\partial x^2}  
      \label{eq17}
\end{equation}

Where  $\pi^x$, $\pi^t$, and $\pi^{xx}$ are extended Infinitesimals, and we use standard procedure to formulate them \cite{cole1951quasi}. These extended infinitesimals can be calculated as:

\begin{equation}
 \pi^x(x, t, u, u_x, u_t) = \eta_x + \eta_u u_x - (\xi_{1})_x u_x - (\xi_{1})_u u^2_x  - (\xi_{2})_x u_t - (\xi_{2})_u  u_x u_t 
  \label{eq18}
\end{equation}

\begin{equation}
  \pi^t(x, t, u, u_x, u_t) = \eta_t + \eta_u u_t - (\xi_{1})_t u_x - (\xi_{2})_u u^2_t  - (\xi_{2})_t u_t - (\xi_{1})_u u_x u_t    
   \label{eq19}
\end{equation}

\begin{equation}
\begin{aligned}
  \pi^{xx}(x, t, u, u_x, u_t, u_{xx}, u_{xt}, u_{tt}) &= \eta_{xx} + [(2\eta_{xu} - (\xi_{1})_{xx})]u_x - (\xi_{2})_{xx} u_t + [{\eta_{uu} - 2 (\xi_{1})_xu}] u_x^2 \\
  &\quad - 2 (\xi_{2})_{xu} u_t u_x - (\xi_{1})_{uu} u_x^3 - (\xi_{2})_{uu} u_t^2 - 2 [(\xi_{2})_x \\
  &\quad + (\xi_{2})_u u_x] u_{xt}+ [(\eta_u - 2 (\xi_{1})_x) - 3 (\xi_{1})_u u_x \\
  &\quad - (\xi_{2})_u u_t u_{xx}]
   \label{eq20}
\end{aligned}
\end{equation}

Suppose if $\mathcal{R} = \frac{\partial u}{\partial t} + a u \frac{\partial u}{\partial x} - v \frac{\partial^2 u}{\partial x^2}$ , then as stated by the infinitesimal criterion of invariance \cite{olver1993applications}, $\mathcal{R}$ remains invariant under the Lie group $G$ if and only if the prolongation of the Lie derivative $\mathbf{L}[\mathcal{R}]$ is zero whenever $\mathcal{R}=0$, for all infinitesimal generators $\mathbf{L}$ of $G$. Mathematically it is expressed as:

\begin{equation}
  \mathbf{L}^{(2)} \mathbf{L}(\mathcal{R}) = 0  \quad
\text{iff} \quad \mathcal{R} = 0
 \label{eq21}
\end{equation}

Thus, the vector field $\mathbf{L}$ along second prolongation $\mathbf{{L}^{(2)}}$ is computed by using eqn \eqref{eq16} and eqn \eqref{eq17} \cite{turgut2017similarity, olver1993applications, li2023utilizing, zhang2023enforcing, freire2010note}:

\begin{equation}
\mathbf{L}^{(2)} \mathbf{L} = \mathbf{L} + \pi^x \frac{d}{du_x} + \pi^t \frac{d}{du_t} + \pi^{xx} \frac{d}{du_{xx}} 
 \label{eq22}
\end{equation}

Furthermore, the coefficients $\eta(x, t, u)$, $\xi_1(x, t, u)$, and $\xi_2(x, t, u)$ are determined using equations \eqref{eq21} in \eqref{eq22}. The detailed computations used to derive these coefficients are omitted here for the sake of clarity and conciseness. For the comprehensive explanation, consult the referenced publications \cite{turgut2017similarity, olver1993applications, li2023utilizing, zhang2023enforcing, freire2010note}. With these coefficients, we can obtain the infinitesimal generators of eqn \eqref{eq1} denoted as $\mathbf{L_1}, \mathbf{L_2}, \ldots, \mathbf{L_5}$. These generators represent the associated vector fields for the one-parameter Lie group of infinitesimal transformations, and they are expressed as:
\begin{equation}
  \mathbf{L_1} = \frac{\partial}{\partial x}
   \label{eq23}
\end{equation}

\begin{equation}
 \mathbf{L_2} = \frac{\partial}{\partial t} 
  \label{eq24}
\end{equation}

\begin{equation}
\mathbf{L_3} = t\frac{\partial}{\partial x} + \frac{\partial}{\partial u}
 \label{eq25}
\end{equation}

\begin{equation}
\mathbf{L_4} = -x\frac{\partial}{\partial x} - 2t\frac{\partial}{\partial t} + u\frac{\partial}{\partial u}  
 \label{eq26}
\end{equation}

\begin{equation}
  \mathbf{L_5} = -tx\frac{\partial}{\partial x} - t^2\frac{\partial}{\partial t} + (tu - x)\frac{\partial}{\partial u} 
   \label{eq27}
\end{equation}
 
These generators, indicate the directions of infinitesimal transformations inside the group and correspond to tangent vectors at the group's identity element. Infinitesimal generators are used to reduce the number of independent variables of PDEs when Lie group action is applied; this process is referred as symmetry reduction. Reducing number of variables often makes it easy to solve PDEs. However, It is important to note that symmetry reduction  can often be a very complex task depending on behaviour of problem. So, Instead of reducing Burgers equation directly, we utilize information provided by infinitesimal generator for Burger's equation along with PINN in novel manner.

 we utilize the Infinitesimal generators of Burger's equation \eqref{eq1} along with  PINN by modifying the loss function, which  differs from previous work \cite{zhang2023enforcing}. Specifically, Instead of using collocation points directly in residual term as defined in \eqref{eq9}, we first transform these points by taking leverage of Infinitesimal generators. Then, another residual term is added to the loss function defined in \eqref{eq32} by using these transformed collocation points. Modified loss function is then optimised to obtained approximate solution. Let's delve deeper for further understanding, which basically begins with the definition of Invariant function.\\
\textbf{Invariant function}
\\
This invariance is frequently exploited to produce solutions to nonlinear equations through construction. Consider a Lie group as defined in \eqref{eq13}, \eqref{eq14} and \eqref{eq15} that transforms \eqref{eq10},

\begin{equation}
\mathcal{R}[x, t, u, u_x, u_t, u_{xx}, u_{xt}, u_{tt}, \ldots]= 0   
 \label{eq28}
\end{equation}

to itself, That is, the variables that have tildes satisfy

\begin{equation}
\mathcal{R}[\tilde{x}, \tilde{t}, \tilde{u}, \tilde{u}_x, \tilde{u}_t, \tilde{u}_{xx}, \tilde{u}_{xt}, \tilde{u}_{tt}, \ldots]= 0
 \label{eq29}
\end{equation}

It follows that if \( u(x, t) \) is a solution of eqn \eqref{eq28}, then \( \tilde{u}(\,\widetilde{x},\,\widetilde{t}\,) \) constructed from the Infinitesimal transformation of Lie group is also a solution of eqn \eqref{eq29}. This is true regardless of whether the equation is linear or nonlinear. Here, we can observe that the differential operator $\mathcal{R}$ maps points on surface to other points on
the same surface \cite{cantwell2002introduction}.
\\
 By employing the PINN, as explained in section 2.1, the loss of residual $\mathcal{R}(t_j^c, x_j^c, \Theta)$ term is minimized by training it on collocation points $\{(t_j^c, x_j^c)\}_{j=1}^{N_c}$, indicating that $\{u(t_j^c, x_j^c)\}_{j=1}^{N_c}$ is a solution of $\mathcal{R}(t_j^c, x_j^c, \Theta)$. Instead of minimising  \(\mathcal{L}_{\text{res}}(\Theta)\)  directly using collocation points as described in eqn \eqref{eq9}, we first transformed these collocation points by using Infinitesimal generator $\mathbf{L}$ of eqn \eqref{eq1}. Mathematically, It can be expressed as
\begin{equation}
    \tilde{x^c} = x^c + \varepsilon \xi_1(t^c,x^c,u(t^c,x^c)) + O(\varepsilon^2) 
     \label{eq30}
\end{equation}

\begin{equation}
     \tilde{t^c} = t^c + \varepsilon \xi_2(t^c,x^c,u(t^c,x^c)) + O(\varepsilon^2)  
      \label{eq31}
\end{equation}

By utilizing these transformed collocation points, i.e., $\{\,(\tilde{t}_j^c, \tilde{x}_j^c)\,\}_{j=1}^{N_c}$, we modify the loss term by adding an extra residual term denoted as \( L_{\text{symm}} \). However, It is important to note that we use same differential operator $\mathfrak{R}$ in \( L_{\text{symm}} \). Mathematically, \( L_{\text{symm}} \) can be expressed as

 \begin{equation}
    \mathcal{L}_{\text{symm}}(\Theta) = \frac{1}{N_c} \sum_{j=1}^{N_c} \left( \mathcal{R}(\tilde{t}_j^c, \tilde{x}_j^c, \Theta) \right)^2 
     \label{eq32}
\end{equation} 

We can obtain modified loss function by incorporating \( L_{\text{symm}} \), represented as

\begin{equation}
\mathcal{L}(\Theta) = \alpha \mathcal{L}_{\text{init}}(\Theta) + \beta \mathcal{L}_{\text{bound}}(\Theta) + \gamma \mathcal{L}_{\text{res}}(\Theta) + \zeta\mathcal{L}_{\text{symm}}(\Theta)  \label{eq33}   
\end{equation}

where $\zeta$ is weight parameter. The incorporation of Lie symmetry into the PINN loss function gives rise to the m-SPINN. It is crucial to remember that we can employ \( n \) \( \mathcal{L}_{\text{symm}} \) terms if we have \( n \) infinitesimal generators $\mathbf{L_1}, \mathbf{L_2}, \ldots, \mathbf{L_n}$. It is essential to carefully choose the infinitesimal generators that are more influential than others or best suited for our particular situation \cite{freire2010note}. Later on, in the results and discussion sections, we will elaborate the efficiency of this technique, examining its performance with Adam optimizer.

\subsection{Modified Adaptive SPINN (m-ASPINN)}
 The performance of neural networks is significantly affected by the activation function selection. By determining whether a neuron should be fired or not, activation functions enable the network to change the input data in nonlinear ways. The network's capacity to tackle complicated issues would be restricted if activation functions were absent as it could only carry out linear changes. Furthermore, activation functions are important to back-propagation because they allow gradients to be evaluated, which is required to update the weights and biases during training. For this reason, choosing the right activation function is crucial to preventing problems like disappearing or expanding gradients. \\
Additionally, the network's ability to manage complicated tasks is also influenced by its size and architecture. Deep learning networks are essential for tackling complex issues, but they may be difficult to train. We need a robust and effective neural network since the loss function becomes increasingly complex as additional terms are added. In order to enhance efficiency of m-SPINN, a recently  developed technique is adopted \cite{jagtap2020adaptive} to modify the activation function by adding the hyper-parameter \(\boldsymbol{\alpha}\) as:

\begin{equation}
    \sigma(\boldsymbol{\alpha} L_k(z^{k-1}))
     \label{eq34}
\end{equation}

 The complete optimization problem leads to finding the minimum of a loss function by optimizing \( \boldsymbol{\alpha} \) along with the weights and biases, i.e., we try to find

\begin{equation}
    \boldsymbol{\alpha}^* = \arg \min_{\alpha \in \mathbb{R}^+ \setminus \{0\}} \mathcal{L}(\boldsymbol{\alpha})
     \label{eq35}
\end{equation}

\( \boldsymbol{\alpha} \) parameter can be updated as

\begin{equation}
    \boldsymbol{\alpha}^{m+1} = \boldsymbol{\alpha}^{m} - \eta \nabla_{\alpha}^{m} \mathcal{L}^{m}(\boldsymbol{\alpha})
     \label{eq36}
\end{equation}

To accelerate convergence, we introduce a constant \( n \), such that \( n \geq 1 \). So the final form of the activation function is given by:

\begin{equation}
    \sigma(\boldsymbol{n\alpha} L_k(z^{k-1}))
     \label{eq37}
\end{equation}

By making activation function adaptive using this technique, m-SPINN is called m-ASPINN. 
\section{Numerical Experiment}
\subsection*{Problem 1}

Consider an example of second order Burger's equation, Since numerous other well-known numerical methods have been used to solve this example, it is convenient for us to evaluate the accuracy of our model and assess the development of PINN model \cite{arora2013numerical}{}.

\begin{equation}
\frac{\partial u(x,t)}{\partial t} = \nu \frac{\partial^2 u(x,t)}{\partial x^2} - u(x,t) \frac{\partial u(x,t)}{\partial x}, \quad \text{where} \quad x \in [0, 1], \quad t \in [0, 3], 
 \label{eq38}
\end{equation}

\begin{equation}
  u(x, 0) = 4x(1-x)  
   \label{eq39}
\end{equation}

\begin{equation}
    u(0, t) = 0 \quad \text{and} \quad u(1, t) = 0.
     \label{eq40}
\end{equation}

Where viscosity, \( v = 0.1 \). The series solution in given for \eqref{eq38}:

\begin{equation}
u(x,t) = \frac{2 \pi v \sum_{n=1}^{\infty} K_n \exp(-n^2 \pi^2 vt) n \sin(n \pi x)}{K_0 + \sum_{n=1}^{\infty} K_n \exp(-n^2 \pi^2 vt) \cos(n \pi x)}
 \label{eq41}
\end{equation}

where:

\begin{equation}
K_0 = \int_{0}^{1} \exp \left\{ -\frac{1}{3v} [x^2 (3 - 2x)] \right\} dx
 \label{eq42}
\end{equation}

\begin{equation}
K_n = 2 \int_{0}^{1} \exp \left\{ -\frac{1}{3v} [x^2 (3 - 2x)] \right\} \cos(n \pi x) dx
 \label{eq43}
\end{equation}

\subsection*{Data}
In this example, we are going to examine three cases: A, B and  C in order to use PINN. Same structural components are used in each case for comparison. Components are given in Table \ref{tab:fixed_components}.

\begin{table}[htbp]
  \centering
  \caption{Fixed Components of the Neural Network}
    \begin{tabular}{ll}
    \toprule
    \textbf{Component} & \textbf{Value} \\
    \midrule
    Activation function & Tanh, GELU , Mish, Swish \\
    Number of layers & 8 \\
    Number of neurons per layer & 40 \\
    Initial condition points, $N_0$ & 500 \\
    Boundary condition points, $N_b$ & 500\\
    Collocation points, $N_r$ & 20,000 \\
    Learning rate schedule & Piecewise Decay \\
    Optimizer & Adam  \\
    Iterations & 50,000  \\
    \bottomrule
    \end{tabular}%
  \label{tab:fixed_components}%
\end{table} 
We use a set of few common data points given in the literature \cite{arora2013numerical},\cite{jiwari2013numerical} and \cite{ozics2003finite}  to evaluate the performance of our PINN model. These points correspond to solutions found using a variety of numerical techniques.\\
Notably, computations are performed using the TensorFlow package. Fixing the seed zero ensures the consistent results for both data points and parameters. Initial, boundary and collocation data points. The 'glorot normal' method is also used to initialize the weights and biases. Furthermore, models are iteratively trained ten times in each case and select the best-performing result among each of them. Later, upon completion of training, approximated solutions are computed. In each case, experiments are conducted by employing various contemporary activation functions, such as Mish \cite{misra1908mish}, Swish\cite{ramachandran2017searching}, hyperbolic tangent (Tanh), exponential hyperbolic tangent (Tanh(exp))\cite{liu2021tanhexp}, Softplus, and Gaussian error linear units (GELU) \cite{hendrycks2016gaussian}. Only Best four performing activation function are selected for analysis. Absolute error is used to make a comparison with numerical techniques Absolute is defined as

\begin{equation}
\text{Absolute Error} = \left| x_{\text{Approximated}} - x_{\text{Exact}} \right|
 \label{eq44}
\end{equation}

\subsection{Results}
\subsubsection*{\textbf{Case A:}}

In this section, solutions of eqn \eqref{eq1} are computed by using conventional PINN model. Residual term can be represented as:

\begin{equation}
\mathcal{L}_{\text{res}}(\Theta) = \sum_{i=1}^{N_r} \left( \frac{\partial u_i}{\partial t_i} + u_i \frac{\partial u_i}{\partial x_i} - v \frac{\partial^2 u_i}{\partial x_i^2} \right)^2
 \label{eq45}
\end{equation}
where $\alpha$=$\beta$=$\gamma$=1.
 Approximated solution by PINN model is given in table   \ref{tab:absolute_error}. This table also  represents the absolute error in between PINN with different activation functions and exact solution. Here  It can be observe that the approximated values from the constructed model of PINN significantly deviate from the exact values and those obtained by other numerical methods \cite{kutluay2004numerical},\cite{arora2013numerical},\cite{jiwari2013numerical} and \cite{ozics2003finite}.
 
\begin{table}[htbp]
\centering
\caption{Comparison of Exact Values and approximated values  of \textbf{PINN} with their Absolute Errors, Solutions are approximated by using different activation functions. }
\scriptsize
\begin{tabular}{|c|c|c|cc|cc|cc|cc|}
\toprule
$x$ & $t$ & \textbf{Exact} & \shortstack{PINN\\(GELU)} & \shortstack{Error\\ (GELU)} & \shortstack{PINN\\(tanh)} & \shortstack{Error \\ (tanh)} &\shortstack{ PINN \\(mish)} & \shortstack{ Error \\ (mish)} & \shortstack{ PINN \\ (swish)} & \shortstack{Error \\ (swish)} \\
\midrule
\multirow{4}{*}{0.25} & 0.4 & 0.31752 & -1.23083 & 1.54835 & -1.16108 & 1.47860 & -2.65378 & 2.97130 & -4.87638 & 5.19390 \\
 & 0.8 & 0.19956 & -0.55792 & 0.75748 & -0.38284 & 0.58240 & -0.47294 & 0.67250 & -0.47342 & 0.67298 \\
 & 1.0 & 0.16560 & -0.39657 & 0.56217 & -0.31683 & 0.48243 & -0.33670 & 0.50230 & -0.32980 & 0.49540 \\
 & 3.0 & 0.02775 & -7.16731 & 7.19506 & -0.12812 & 0.15587 & -0.40639 & 0.43414 & 0.21507 & 0.18732 \\
\midrule
\multirow{4}{*}{0.50} & 0.4 & 0.58454 & 0.39931 & 0.18523 & 0.24695 & 0.33759 & 0.68999 & 0.10546 & 0.78829 & 0.20375 \\
 & 0.8 & 0.36740 & 0.00008 & 0.36732 & 0.00020 & 0.36720 & 0.00003 & 0.36737 & 0.00006 & 0.36734 \\
 & 1.0 & 0.29834 & 0.00069 & 0.29765 & -0.00051 & 0.29885 & -0.00005 & 0.29839 & -0.00009 & 0.29843 \\
 & 3.0 & 0.04106 & 7.69596 & 7.65490 & -0.02826 & 0.06932 & -0.31882 & 0.35988 & 0.09757 & 0.05651 \\
\midrule
\multirow{4}{*}{0.75} & 0.4 & 0.64562 & 1.07120 & 0.42558 & 1.02510 & 0.37948 & 0.99708 & 0.35146 & 0.96786 & 0.32224 \\
 & 0.8 & 0.38534 & 0.46166 & 0.07632 & 0.45796 & 0.07262 & 0.45798 & 0.07264 & 0.45798 & 0.07264 \\
 & 1.0 & 0.29586 & 0.30290 & 0.00704 & 0.29834 & 0.00248 & 0.29834 & 0.00249 & 0.29835 & 0.00249 \\
 & 3.0 & 0.03044 & 40.53058 & 40.50014 & 0.04203 & 0.01159 & -0.44772 & 0.47816 & -0.32683 & 0.35727 \\
\bottomrule
\end{tabular}
\label{tab:absolute_error}
\end{table}
 Due to large discrepancy, improvements are made in PINN model that are discussed in case B.
\subsubsection*{\textbf{Case B:}}
In this case, we modify case A by making changes to the loss function. Instead of using 20,000 collocation points (\(t_r\), \(x_r\)) directly for the residual term, Firstly  these points are transformed using infinitesimal transformations \eqref{eq13} and \eqref{eq14}.  (\(\tilde{t}_r\), \(\tilde{x}_r\)) are transformed points. The values of \(\xi_1\) and \(\xi_2\)  are obtained using equation \eqref{eq27}, where \(\varepsilon = 0.5\). The transformed collocation points are then used to evaluate \(L_{\text{symm}}\), as defined in equation \eqref{eq32}. \(L_{\text{symm}}\) is subsequently include in the loss function of case A. Finally  loss function is minimised through training. Here, in this case, by employing Lie symmetry, we refer to our PINN model as the  m-SPINN.
Table \ref{tab:mspinn_comparison} represents the approximated solution by trained m-SPINN model. It is clear from Table \ref{tab:mspinn_comparison} that m-SPINN is far more better than PINN. The performance of all m-SPINN across all activation functions produced approximately accurate results  and close to exact and other numerical techniques. Due to superior performance of GELU comparatively to other activation functions during experiment, GELU is proceeded in next section to enhance our results more. Adaptive activation function technique is utilized for further enhancement of results which is explained in next case (case C).

\begin{table}[htbp]
\centering
\caption{Comparison of Exact Values and approximated values  of \textbf{m-SPINN} with their Absolute Errors, Solutions are approximated by using different activation functions }
\scriptsize 
\begin{tabular}{|c|c|c|cc|cc|cc|cc|}
\toprule
$x$ & $t$ & \textbf{Exact} & \shortstack{m-SPINN\\(GELU)} & \shortstack{Error\\ (GELU)} & \shortstack{m-SPINN\\(tanh)} & \shortstack{Error \\ (tanh)} &\shortstack{ m-SPINN \\(mish)} & \shortstack{ Error \\ (mish)} & \shortstack{ m-SPINN \\ (swish)} & \shortstack{Error \\ (swish)} \\
\midrule
\multirow{4}{*}{0.25} & 0.4 & 0.31752 & 0.317410 & 0.000110 & 0.317062 & 0.000458 & 0.317612 & 0.000092 & 0.317399 & 0.000121 \\
 & 0.8 & 0.19956 & 0.199625 & 0.000065 & 0.199678 & 0.000118 & 0.199697 & 0.000137 & 0.199769 & 0.000209 \\
 & 1.0 & 0.16560 & 0.165572 & 0.000028 & 0.165735 & 0.000135 & 0.165631 & 0.000031 & 0.165656 & 0.000056 \\
 & 3.0 & 0.02775 & 0.027663 & 0.000087 & 0.027659 & 0.000091 & 0.027799 & 0.000049 & 0.027600 & 0.000150 \\
\midrule
\multirow{4}{*}{0.50} & 0.4 & 0.58454 & 0.584480 & 0.000060 & 0.584365 & 0.000175 & 0.584405 & 0.000135 & 0.584567 & 0.000027 \\
 & 0.8 & 0.36740 & 0.367376 & 0.000024 & 0.367263 & 0.000137 & 0.367538 & 0.000138 & 0.367457 & 0.000057 \\
 & 1.0 & 0.29834 & 0.298321 & 0.000019 & 0.298372 & 0.000032 & 0.298485 & 0.000145 & 0.298379 & 0.000039 \\
 & 3.0 & 0.04106 & 0.041110 & 0.000050 & 0.041023 & 0.000037 & 0.041252 & 0.000192 & 0.041178 & 0.000118 \\
\midrule
\multirow{4}{*}{0.75} & 0.4 & 0.64562 & 0.645554 & 0.000066 & 0.645517 & 0.000103 & 0.645511 & 0.000109 & 0.645577 & 0.000043 \\
 & 0.8 & 0.38534 & 0.385231 & 0.000109 & 0.385130 & 0.000210 & 0.385442 & 0.000102 & 0.385374 & 0.000034 \\
 & 1.0 & 0.29586 & 0.295817 & 0.000043 & 0.295860 & 0.000000 & 0.296039 & 0.000179 & 0.295804 & 0.000056 \\
 & 3.0 & 0.03044 & 0.030449 & 0.000009 & 0.030488 & 0.000048 & 0.030587 & 0.000147 & 0.030518 & 0.000078 \\
\bottomrule
\end{tabular}
\label{tab:mspinn_comparison}
\end{table}

\subsubsection*{\textbf{Case C:}}

\begin{figure}[htbp] 
    \centering
    \caption{Absolute Error Comparison between m-ASPINN, WA-DQM,and MCB-DQM}
    \includegraphics[width=1
    \textwidth]{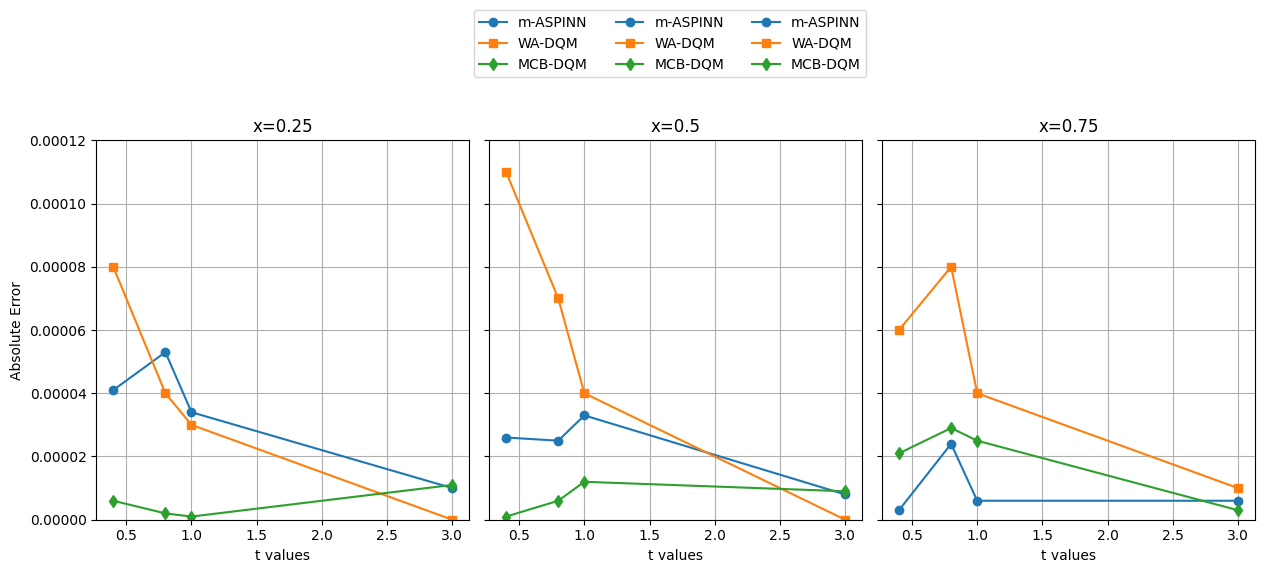} %
    \label{aspinn}
    \end{figure}

\begin{figure}
    \centering
    \caption{ Physical behaviour of m-ASPINN approximated solution of u(x,t) at v = 0.1 (right) and series solution (left).}
    \begin{subfigure}[b]{0.45\textwidth}
        \centering
        \includegraphics[width=\textwidth]{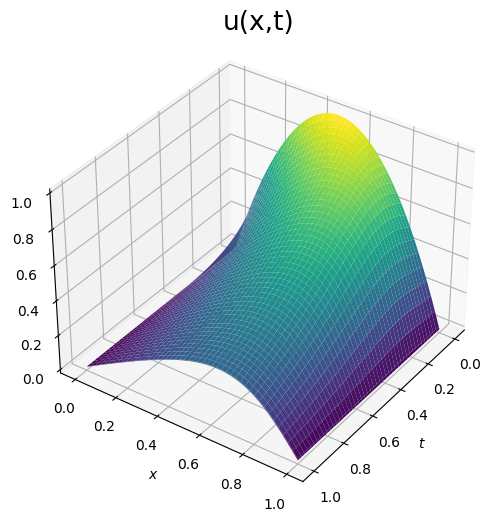} 
        \caption{Series solution surface}
        \label{ac}
    \end{subfigure}
    \hfill
    \begin{subfigure}[b]{0.45\textwidth}
        \centering
        \includegraphics[width=\textwidth]{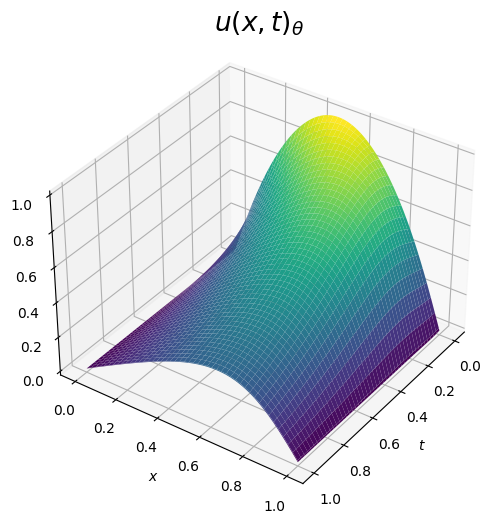} 
        \caption{m-ASPINN solution surface}
        \label{cross curve}
    \end{subfigure}

    \label{ajk}
\end{figure}

\begin{figure}
    \centering
    \caption{ Cross-section curves by m-ASPINN approximated solution of u(x,t) at v = 0.1 (right) at t= 0, 0.4, 0.8 and 1  and  series solution approximated cross-section curves (left). }
    \begin{subfigure}[b]{0.45\textwidth}
        \centering
        \includegraphics[width=\textwidth]{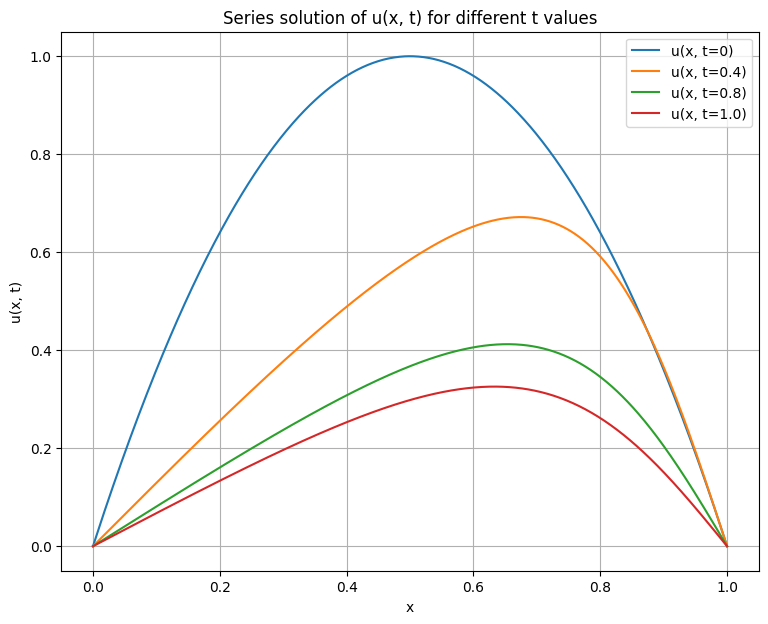} 
        \caption{Series solution cross-section curves}
        \label{ac1}
    \end{subfigure}
    \hfill
    \begin{subfigure}[b]{0.45\textwidth}
        \centering
        \includegraphics[width=\textwidth]{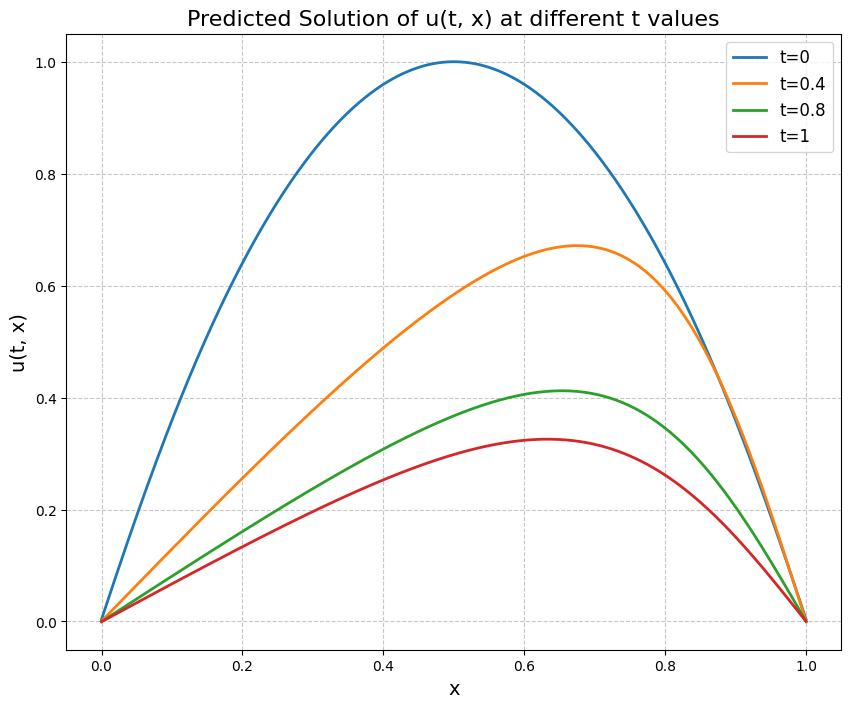} 
        \caption{m-ASPINN cross-section curves}
        \label{cross curve1}
    \end{subfigure}

    \label{ajk1}
\end{figure}

\begin{figure}
    \centering
    \caption{ Gradient lines plot by m-ASPINN method (right) and by series solution (left) for solution surfaces as shown in figure \ref{ajk}}
    \begin{subfigure}[b]{0.45\textwidth}
        \centering
        \includegraphics[width=\textwidth]{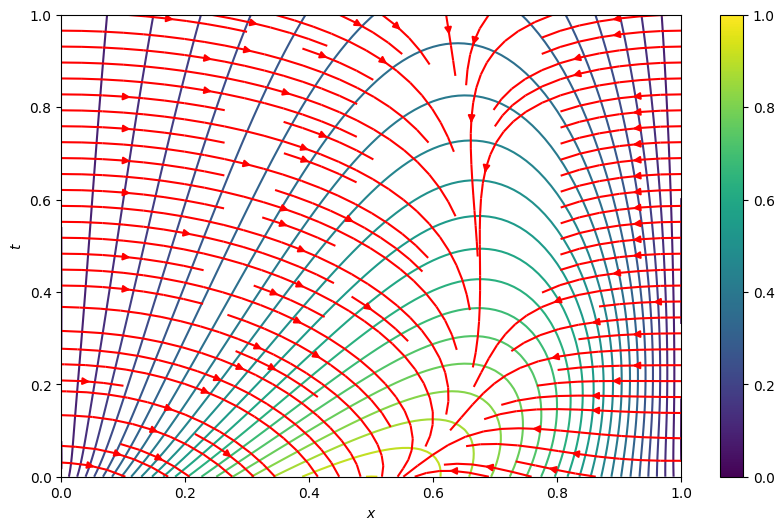} 
        
        \label{heatmap}
    \end{subfigure}
    \hfill
    \begin{subfigure}[b]{0.45\textwidth}
        \centering
        \includegraphics[width=\textwidth]{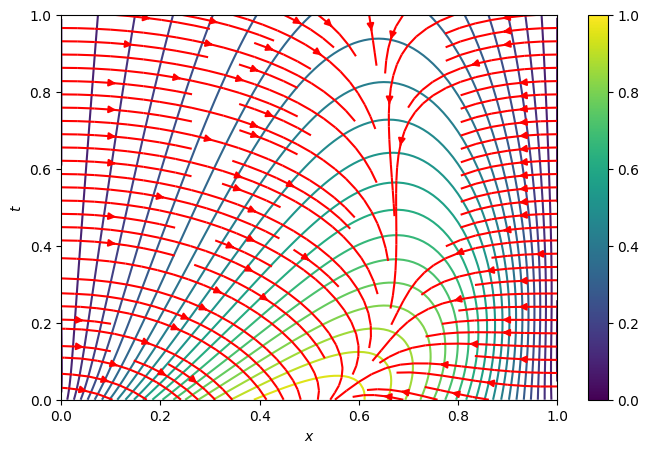} 
    
        \label{Contour lines}
    \end{subfigure}
    \label{Heatcontour}
\end{figure}

\begin{figure}[htbp]
    \centering
    \begin{minipage}{0.45\textwidth}
        \centering
        \caption{Number of data points and corresponding computation times for problem 1.}
        \begin{tabular}{|c|c|}
        \hline
    \textbf{\shortstack{Number of \\ Points}} & \textbf{\shortstack{ Computation Time \\ (seconds)}} \\
        \hline
        1000 & 8.498138 \\
        \hline
        2000 & 14.448773 \\
        \hline
        3000 & 21.888418 \\
        \hline
        4000 & 28.745561 \\
        \hline
        5000 & 36.741615 \\
        \hline
        6000 & 43.736098 \\
        \hline
        7000 & 50.661646 \\
        \hline
        8000 & 58.080971 \\
        \hline
        9000 & 66.747062 \\
        \hline
        10000 & 71.743773 \\
        \hline
        \end{tabular}
\label{table:computation_time 1}
    \end{minipage}
    \hfill
    \begin{minipage}{0.45\textwidth}
        \centering
        \caption{Plot of computation time vs. number of data points.} \includegraphics[width=\textwidth]{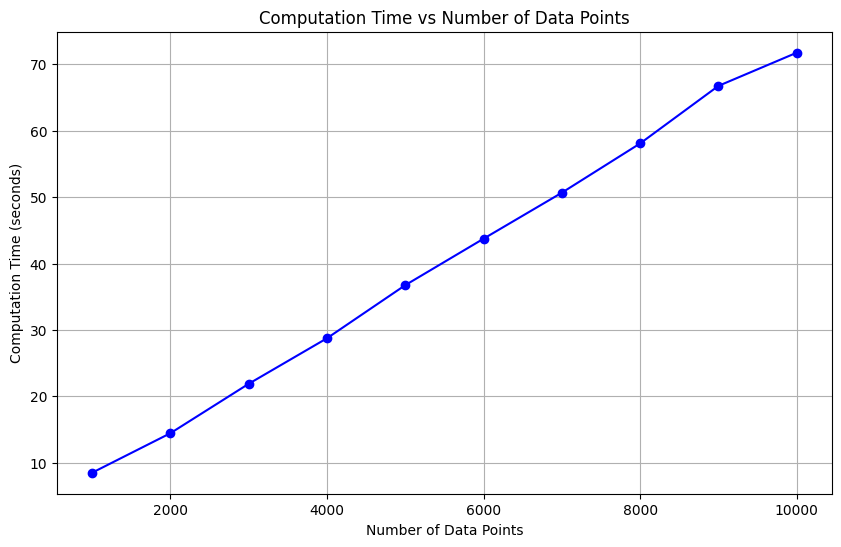}
    \label{fig:computation_time_plot 2}
    \end{minipage}
\end{figure}

\begin{figure}[htbp]
	\centering
	\captionsetup{justification=centering}
	\caption{Solution surface behaviour in two dimensional heat map by m-ASPINN (right) and by series solution (left)}
	\begin{subfigure}[b]{0.45\textwidth}
		\centering
		\includegraphics[width=\textwidth]{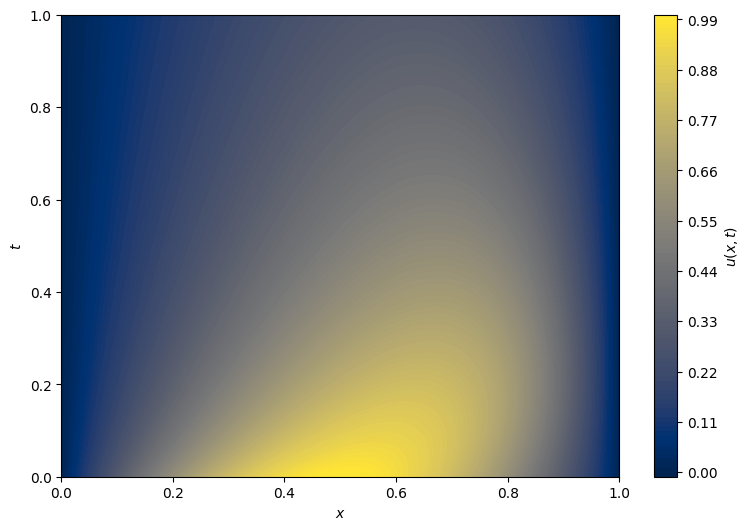}
		\caption{Series solution}
		\label{heatmap2}
	\end{subfigure}
	\hfill
	\begin{subfigure}[b]{0.47\textwidth}
		\centering
		\includegraphics[width=\textwidth]{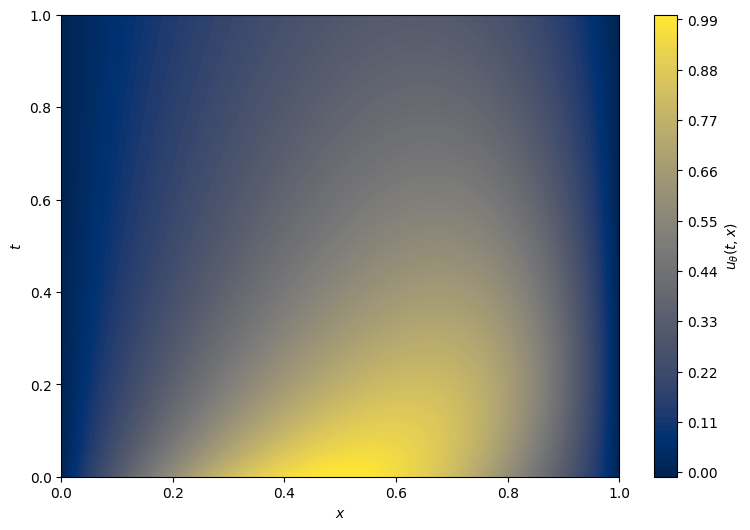}
		\caption{m-ASPINN}
		\label{Contour_lines12}
	\end{subfigure}
	\label{Heatcontour1}
\end{figure}

To further improve the efficiency of our results in m-SPINN, we employ adaptive activation function  as defined in subsections 2.3. A hyper parameter \(\boldsymbol{\alpha}\)  = 0.1 is used in this technique \cite{jagtap2020adaptive}.
The modified SPINN model, now called modified adaptive SPINN or m-ASPINN. Table \ref{tab:tab_activation} shows the m-ASPINN results. Table \ref{tab:solutions_comparison} shows the approximate solutions of m-ASPINN and other numerical methods. Absolute errors in Table \ref{vv} clearly demonstrate that m-ASPINN is more accurate than the m-SPINN model explained in case B and it is almost as efficient as other numerical methods, specifically MCB-DQM \cite{arora2013numerical}.
Figure \ref{aspinn} provides graphical illustration of Tables \ref{vv} . We eliminate LS-QB-SPLINE-FEM in Figure \ref{aspinn} because of its slightly larger deviation as compared other to other techniques in order to visualise more closely. Figure \ref{ajk} represents the 3D surface solution of \eqref{eq1} by m-ASPINN method and by series solution. 
Figure \ref{ajk1} shows the cross-section curves of solution surface for different values of t, cross curves are obtained for t = 0, 0.4, 0.8 , 1 by using m-ASPINN and Series solution method.
Figure \ref{Heatcontour} shows the gradient line plot that illustrates the behavior of the surface
u with respect to the spatial dimension and time. While, Figure \ref{Heatcontour1} shows the heat map of  u in 2D plot, which is a nice way to understand that how solution u evolves with time. After training m-ASPINN model, the computation time of model is evaluated by generating data points ranging from 1,000 to 10,000 in intervals of 1,000. Figure \ref{fig:computation_time_plot 2} shows the computational efficiency of the m-ASPINN model.

\begin{table}[htbp]
	\centering
	\renewcommand{\arraystretch}{1.2} 
	\caption{Comparison of Exact Values and approximated values of \textbf{m-ASPINN} with their Absolute Errors, Solutions are approximated by using different activation functions}
	\scriptsize
	\begin{tabular}{|c|c|c|cc|cc|cc|cc|}
		\hline
		$x$ & $t$ & \textbf{Exact} & \shortstack{\textbf{m-ASPINN}\\(GELU)} & \shortstack{\textbf{Error}\\ (GELU)} & \shortstack{m-ASPINN\\(tanh)} & \shortstack{Error \\ (tanh)} &\shortstack{m-ASPINN \\(mish)} & \shortstack{ Error \\ (mish)} & \shortstack{ m-ASPINN \\ (swish)} & \shortstack{Error \\ (swish)} \\
		\hline
		\multirow{4}{*}{0.25} & 0.4 & 0.31752 & 0.317561 & 0.000 041 & 0.317103 & 0.000417 & 0.317615 & 0.000089 & 0.317446 & 0.000074  \\
		& 0.8 & 0.19956 & 0.199613 & 0.000 053 & 0.199437 & 0.000123 & 0.199593 & 0.000033 & 0.19955 & 0.000010 \\
		& 1.0 & 0.16560 & 0.165634 & 0.000 034 & 0.165500 & 0.000100 & 0.165629 & 0.000029 & 0.165535 & 0.000065 \\
		& 3.0 & 0.02775 & 0.027663 & 0.000 010 & 0.027654 & 0.000096 & 0.027791 & 0.000041 & 0.02776 & 0.000034 \\
		\hline
		\multirow{4}{*}{0.50} & 0.4 & 0.58454 & 0.584480 & 0.000 026 & 0.584436 & 0.000104 & 0.58492 & 0.000038 & 0.584510 & 0.000030 \\
		& 0.8 & 0.36740 & 0.367425 & 0.000 025 & 0.367261 & 0.000139 & 0.367517 & 0.000117 & 0.367346 & 0.000054 \\
		& 1.0 & 0.29834 & 0.298373& 0.000 033 & 0.298368 & 0.000028 & 0.298301 & 0.000039 & 0.29826 & 0.000080 \\
		& 3.0 & 0.04106 & 0.041110 & 0.000 008 & 0.041011 & 0.000049 & 0.041138 & 0.000078 & 0.041068 & 0.000022 \\
		\hline
		\multirow{4}{*}{0.75} & 0.4 & 0.64562 & 0.645617 & 0.000 003 & 0.645567 & 0.000053 & 0.645522 & 0.000098 & 0.645648 & 0.000028 \\
		& 0.8 & 0.38534 & 0.385364 & 0.000 024 & 0.385282 & 0.000058 & 0.385453 & 0.000113 & 0.385385 & 0.000045 \\
		& 1.0 & 0.29586 & 0.295866 & 0.000 006 & 0.295814 & 0.000046 & 0.296025 & 0.000165 & 0.295704 & 0.000090 \\
		& 3.0 & 0.03044 & 0.030449 & 0.000 006  & 0.030446 & 0.000015 & 0.030563 & 0.000123 & 0.030426 & 0.000014 \\
		\hline
	\end{tabular}
	\label{tab:tab_activation}
\end{table}

\begin{table}[htbp]
  \centering
 
  \caption{Approximated results of m-ASPINN, MCB-DQM\cite{arora2013numerical} WA-DQM\cite{jiwari2013numerical}, LS-QB spline-FEM\cite{kutluay2004numerical} and exact solution, where \(v = 0.1\)}
    \begin{tabular}{|c|c|c|c|c|c|c|c|c|}
    \toprule
    \textbf{x} & \textbf{t} & \textbf{Exact} & \textbf{\shortstack{\cite{arora2013numerical} \\ MCB-\\DQM}} & \textbf{\shortstack{\cite{jiwari2013numerical} \\ WA-\\DQM}} & \textbf{\shortstack{\cite{kutluay2004numerical}\\LS-QB\\Spline-FEM}} & \textbf{m-ASPINN} \\
    \midrule
    0.25 & 0.4 &  0.31752 & 0.317526 &  0.31744 &  0.32091 & 0.317561 \\
     & 0.8 &  0.19956 & 0.199558 &  0.19952 &  0.20211 & 0.199613 \\
     & 1 &  0.16560 &  0.165601 &  0.16557 &  0.16782 & 0.165634 \\
      & 3 &  0.02775 &  0.027761 &  0.02775 &  0.02828 & 0.02776 \\
    \midrule
   0.5 & 0.4 &  0.58454 & 0.584541 &  0.58443 &  0.58788 & 0.584566 \\
     & 0.8 &  0.36740 &  0.367406 &  0.36773 & 0.37111 &  0.367425 \\
     & 1 &  0.29834 &  0.298352 &  0.29830 &  0.30183 & 0.298373 \\
    & 3 &  0.04106 &  0.041069 &  0.04106 &  0.04185 & 0.041068 \\
    \midrule
    0.75 & 0.4 &  0.64562 & 0.645641 & 0.64556 &  0.65054 & 0.645617 \\
  & 0.8 &  0.38534 &  0.385369 &  0.38526 &   0.39068 & 0.385364 \\
     & 1 &  0.29586 &  0.295885 &  0.29582 &  0.30057 & 0.295866 \\
          & 3 &  0.03044 &  0.030443 &  0.03043 &  0.03106  & 0.030446 \\
    \bottomrule
    \end{tabular}%
  \label{tab:solutions_comparison}%
\end{table}

\begin{table}[htbp]
  \centering
  \caption{Point-wise Absolute Error Comparison of m-ASPINN with numerical methods}
  \begin{tabular}{p{2cm} *{4}{S[table-format=1.6]} S[table-format=1.6]}
    \toprule
    \textbf{x} & \textbf{t} & \textbf{MCB-DQM \cite{arora2013numerical}} & \textbf{m-ASPINN} & \textbf{WA-DQM \cite{jiwari2013numerical}} & \textbf{LS-QB Spline-FEM \cite{kutluay2004numerical}} \\
    \midrule
    0.25 & 0.4 & 0.000006 & 0.000041 & 0.000080 & 0.003309 \\
         & 0.8 & 0.000002 & 0.000053 & 0.000040 & 0.002550 \\
         & 1   & 0.000001 & 0.000034 & 0.000030 & 0.002220 \\
         & 3   & 0.000011 & 0.000010 & 0.000000 & 0.000530 \\
    \midrule
    0.5  & 0.4 & 0.000001 & 0.000026 & 0.000110 & 0.003340 \\
         & 0.8 & 0.000006 & 0.000025 & 0.000070 & 0.003710 \\
         & 1   & 0.000012 & 0.000033 & 0.000040 & 0.003490 \\
         & 3   & 0.000009 & 0.000008 & 0.000000 & 0.000790 \\
    \midrule
    0.75 & 0.4 & 0.000021 & 0.000003 & 0.000060 & 0.004920 \\
         & 0.8 & 0.000029 & 0.000024 & 0.000080 & 0.005340 \\
         & 1   & 0.000025 & 0.000006 & 0.000040 & 0.004710 \\
         & 3   & 0.000030 & 0.000006 & 0.000040 & 0.004710 \\
    \bottomrule
  \end{tabular}
  \label{vv}
\end{table}

\begin{figure}[htbp]
	\centering
	\captionsetup{justification=centering}
	\caption{Solution surface behaviour in two dimensional heat map by m-ASPINN (right) and by series solution (left)}
	\begin{subfigure}[b]{0.45\textwidth}
		\centering
		\includegraphics[width=\textwidth]{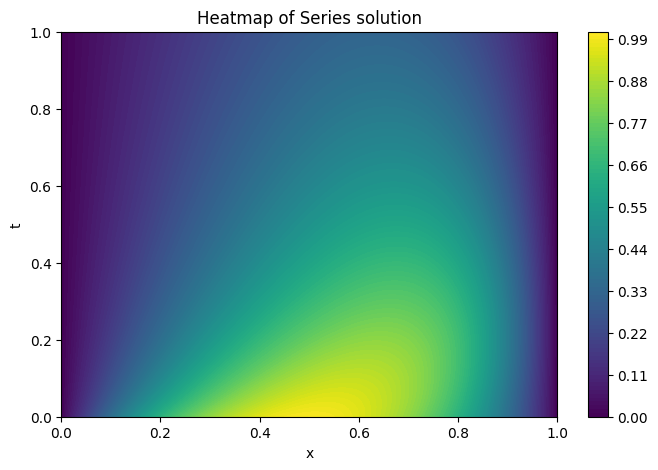}
		\caption{Series solution}
		\label{heatmap1}
	\end{subfigure}
	\hfill
	\begin{subfigure}[b]{0.47\textwidth}
		\centering
		\includegraphics[width=\textwidth]{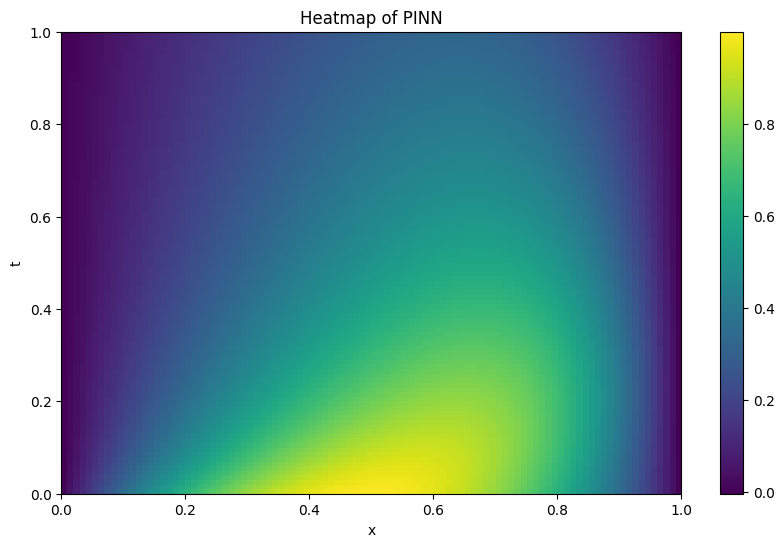}
		\caption{m-ASPINN}
		\label{Contour_lines1}
	\end{subfigure}
    
\end{figure}

\begin{figure}[htbp] 
    \centering
    \caption{Absolute Error surface for problem 1, where $\delta{t}=0.01$, $h= 0.01$ }
    \includegraphics[width=0.7
    \textwidth]{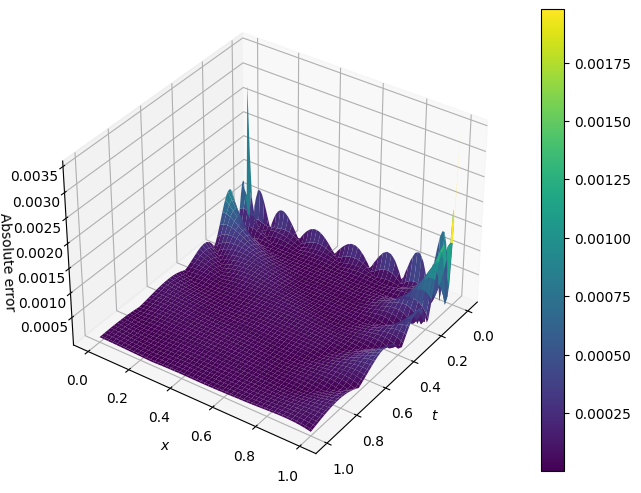} %
    \label{Absolute error surface}
    \end{figure}

\subsection*{Problem 2}

The partial differential equation is given by \cite{mukundan2015efficient}:
\begin{equation}
\frac{\partial u(x,t)}{\partial t} = \nu \frac{\partial^2 u(x,t)}{\partial x^2} - u(x,t) \frac{\partial u(x,t)}{\partial x}
 \label{eq46}
\end{equation}

with \( x \) in the interval \([0,1]\) and \( t \) in the interval \([0,T]\).

Initial condition:
\begin{equation}
u(x,0) = \sin(\pi x), \quad 0 < x < 1
 \label{eq47}
\end{equation}

Boundary conditions:
\begin{equation}
u(0,t) = 0 = u(1,t), \quad 0 \le t \le T
 \label{eq48}
\end{equation}

Series solution (exact):
\begin{equation}
u(x,t) = \frac{2 \pi \nu \sum_{n=1}^{\infty} M_n \exp(-n^2 \pi^2 \nu t) n \sin(n \pi x)}{M_0 + \sum_{n=1}^{\infty} M_n \exp(-n^2 \pi^2 \nu t) \cos(n \pi x)}
 \label{eq49}
\end{equation}
where:
\begin{align}
M_0 &= \int_0^1 \exp \left\{ -\frac{1}{2 \pi \nu} [1 - \cos(\pi x)] \right\} dx \\
M_n &= 2 \int_0^1 \exp \left\{ -\frac{1}{2 \pi \nu} [1 - \cos(\pi x)] \right\} \cos(n \pi x) dx
\end{align}

\subsection*{Data}
In Problem 2, the same procedure and data as in Problem 1 are utilized. Also, in this problem, the Simple PINN model fails to understand the relationship between dependent and independent variables \eqref{eq46}. An infinitesimal transformation of the Lie group is then applied to enhance the results, followed by adaptivity on the GELU activation function, as it is performed in Problem 1. Here, only the results of m-ASPINN (case C) are displayed and compared with the backward differentiation formula (BDF) and its variants, as presented in \cite{mukundan2015efficient}, against the exact solution.

\begin{table}[htbp]
\centering
\caption{Comparison of computed solution approximated by m-ASPINN with BDF-1, BDF-2, BDF-3 and the exact solution at different times and spatial dimension, where BDF\cite{mukundan2015efficient} are discretized at $\nu = 0.1$, $\Delta x = 0.0125$ and $\Delta t = 0.01$.}
\begin{tabular}{ccccccc}
\toprule
$x$ & $T$ & \multicolumn{3}{c}{Computed solution \cite{mukundan2015efficient}} & \multicolumn{1}{c}{Exact solution} & \textbf{m-ASPINN} \\
\cmidrule(r){3-5} \cmidrule(r){6-6}
    &     & BDF-1 & BDF-2 & BDF-3 &  \\
\midrule
0.25 & 2.4 & $4.807 \times 10^{-2}$ & $4.756 \times 10^{-2}$ & $4.755 \times 10^{-2}$ & $4.755 \times 10^{-2}$ & $4.757 \times 10^{-2}$ \\
     & 2.6 & $4.003 \times 10^{-2}$ & $3.956 \times 10^{-2}$ & $3.955 \times 10^{-2}$ & $3.955 \times 10^{-2}$ & $3.956 \times 10^{-2}$ \\
     & 3.0 & $2.759 \times 10^{-2}$ & $2.721 \times 10^{-2}$ & $2.720 \times 10^{-2}$ & $2.720 \times 10^{-2}$ & $2.720 \times 10^{-2}$ \\
0.5  & 2.4 & $7.354 \times 10^{-2}$ & $7.270 \times 10^{-2}$ & $7.268 \times 10^{-2}$ & $7.269 \times 10^{-2}$ & $7.269 \times 10^{-2}$ \\
     & 2.6 & $6.043 \times 10^{-2}$ & $5.968 \times 10^{-2}$ & $5.966 \times 10^{-2}$ & $5.967 \times 10^{-2}$ & $5.967 \times 10^{-2}$ \\
     & 3.0 & $4.080 \times 10^{-2}$ & $4.021 \times 10^{-2}$ & $4.020 \times 10^{-2}$ & $4.020 \times 10^{-2}$ & $4.020 \times 10^{-2}$ \\
0.75 & 2.4 & $5.664 \times 10^{-2}$ & $5.594 \times 10^{-2}$ & $5.593 \times 10^{-2}$ & $5.593 \times 10^{-2}$ & $5.594 \times 10^{-2}$ \\
     & 2.6 & $4.582 \times 10^{-2}$ & $4.522 \times 10^{-2}$ & $4.520 \times 10^{-2}$ & $4.521 \times 10^{-2}$ & $4.521 \times 10^{-2}$ \\
     & 3.0 & $3.023 \times 10^{-2}$ & $2.978 \times 10^{-2}$ & $2.977 \times 10^{-2}$ & $2.977 \times 10^{-2}$ & $2.979 \times 10^{-2}$ \\
\bottomrule
\end{tabular}
\end{table}

\begin{table}[htbp]
\centering
\caption{Comparison of computed solution approximated by m-ASPINN with BDF-1, BDF-2, BDF-3 and the exact solution at different times and spatial dimension, where BDF\cite{mukundan2015efficient} are computed at different space points for example at $T=2.5$ for $\nu = 0.1$ and $\Delta t = 0.01$.}
\begin{tabular}{cccccc}
\toprule
$x$ & \multicolumn{3}{c}{Computed solution \cite{mukundan2015efficient}} & Exact solution & \textbf{m-ASPINN} \\
\cmidrule(r){2-4} \cmidrule(r){5-5}
    & BDF-1 & BDF-2 & BDF-3  \\
\midrule
0.1 & $1.873 \times 10^{-2}$ & $1.852 \times 10^{-2}$ & $1.852 \times 10^{-2}$ & $1.852 \times 10^{-2}$ & $1.852 \times 10^{-2}$ \\
0.2 & $3.611 \times 10^{-2}$ & $3.571 \times 10^{-2}$ & $3.571 \times 10^{-2}$ & $3.571 \times 10^{-2}$ & $3.572 \times 10^{-2}$ \\
0.3 & $5.080 \times 10^{-2}$ & $5.022 \times 10^{-2}$ & $5.021 \times 10^{-2}$ & $5.021 \times 10^{-2}$ & $5.022 \times 10^{-2}$ \\
0.4 & $6.141 \times 10^{-2}$ & $6.070 \times 10^{-2}$ & $6.068 \times 10^{-2}$ & $6.069 \times 10^{-2}$ & $6.069 \times 10^{-2}$ \\
0.5 & $6.666 \times 10^{-2}$ & $6.587 \times 10^{-2}$ & $6.585 \times 10^{-2}$ & $6.586 \times 10^{-2}$ & $6.586 \times 10^{-2}$ \\
0.6 & $6.553 \times 10^{-2}$ & $6.473 \times 10^{-2}$ & $6.471 \times 10^{-2}$ & $6.471 \times 10^{-2}$ & $6.472 \times 10^{-2}$ \\
0.7 & $5.748 \times 10^{-2}$ & $5.676 \times 10^{-2}$ & $5.674 \times 10^{-2}$ & $5.675 \times 10^{-2}$ & $5.676 \times 10^{-2}$ \\
0.8 & $4.283 \times 10^{-2}$ & $4.227 \times 10^{-2}$ & $4.226 \times 10^{-2}$ & $4.226 \times 10^{-2}$ & $4.227 \times 10^{-2}$ \\
0.9 & $2.289 \times 10^{-2}$ & $2.259 \times 10^{-2}$ & $2.258 \times 10^{-2}$ & $2.258 \times 10^{-2}$ & $2.259 \times 10^{-2}$ \\
\bottomrule
\end{tabular}
\end{table}

\begin{figure}
    \centering
    \caption{ Physical behaviour of m-ASPINN approximated solution at v = 0.1 (right) and series solution (left) for problem 2.}
    \begin{subfigure}[b]{0.45\textwidth}
        \centering
        \includegraphics[width=\textwidth]{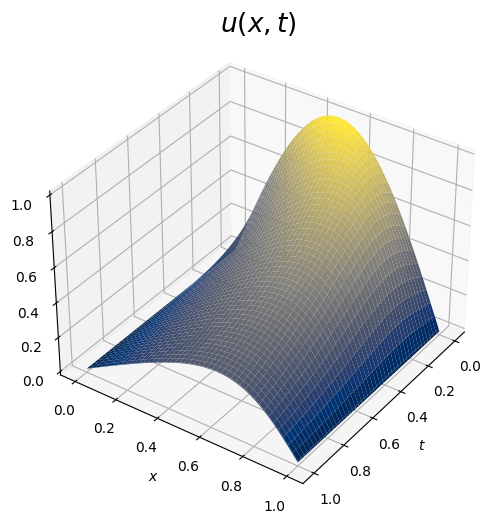} 
        \caption{Series solution surface}
        \label{Series 3D}
    \end{subfigure}
    \hfill
    \begin{subfigure}[b]{0.45\textwidth}
        \centering
        \includegraphics[width=\textwidth]{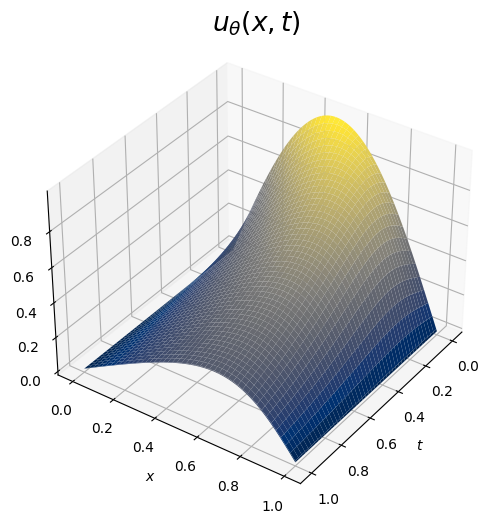} 
        \caption{m-ASPINN solution surface}
        \label{m-ASPINN 3D}
    \end{subfigure}

    \label{3D Surface}
\end{figure}

\begin{figure}
    \centering
    \caption{ Gradient lines plot by m-ASPINN method (right) and by series solution (left) for problem 2 }
    \begin{subfigure}[b]{0.45\textwidth}
        \centering
        \includegraphics[width=\textwidth]{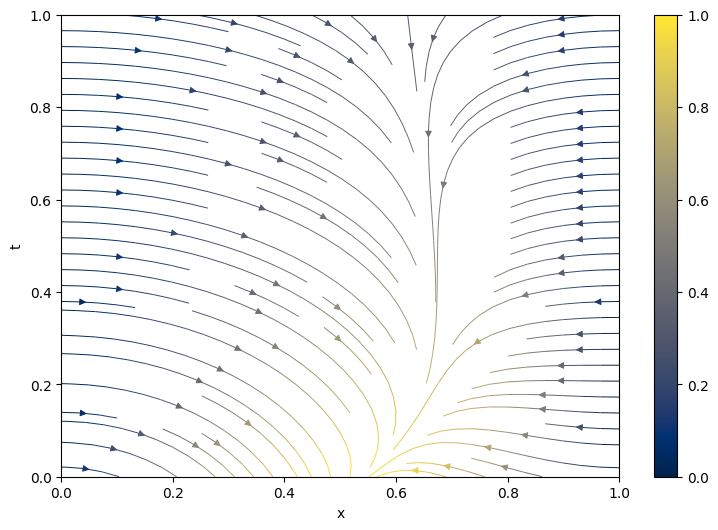} 
        
        \label{series streamlines}
    \end{subfigure}
    \hfill
    \begin{subfigure}[b]{0.45\textwidth}
        \centering
        \includegraphics[width=\textwidth]{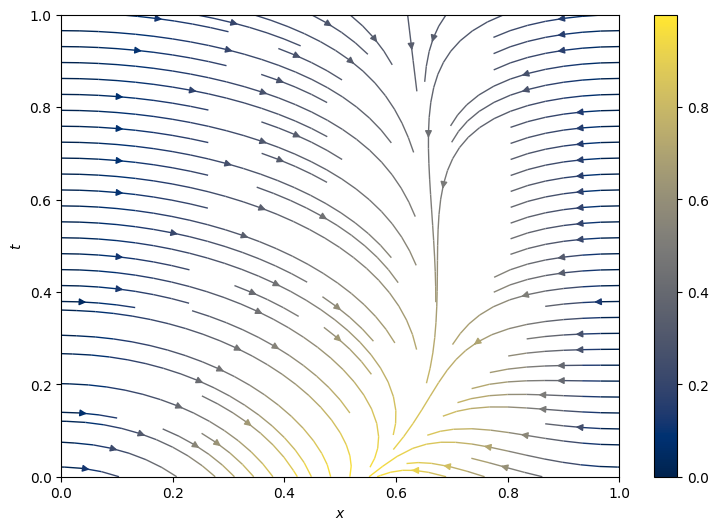} 
    
        \label{m-ASPINN stream lines}
    \end{subfigure}
    \label{stream lines}
\end{figure}

\begin{figure}
    \centering
    \caption{ This shows the comparison of Cross-section curves by m-ASPINN and series solution at t = 1.7, 2.4, 2.5, 2.6 and 3.00 for v= 0.1 (left) and It also illustrates the
Absolute error loss surface for \(\Delta t = 0.01\) and \(\Delta x = 0.01\). where $L_2$ norm is defined in \cite{mukundan2015efficient}. }

    \begin{subfigure}[b]{0.45\textwidth}
        \centering
        \includegraphics[width=\textwidth]{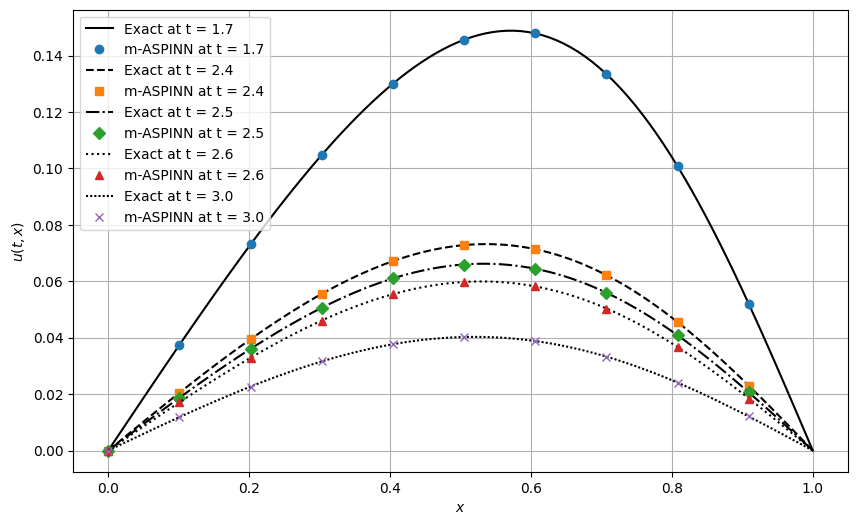} 
        \caption{cross-section curves}
        \label{cross-section curves}
    \end{subfigure}
    \hfill
    \begin{subfigure}[b]{0.45\textwidth}
        \centering
        \includegraphics[width=\textwidth]{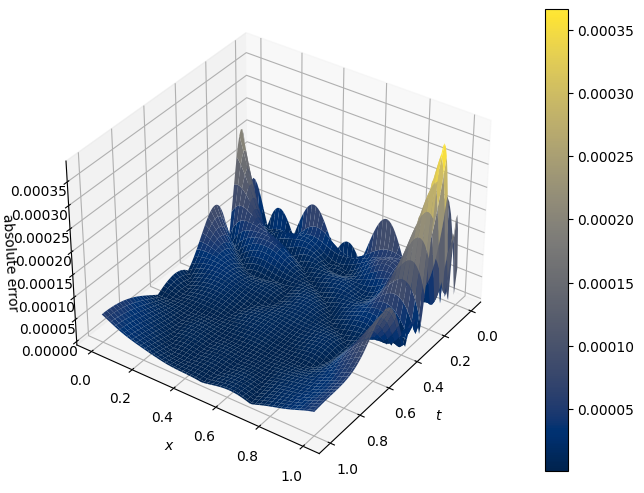} 
        \caption{Absolute error surface}
        \label{ $L_2$ 
   norm loss function}
    \end{subfigure}

    \label{loss and cross}
\end{figure}

\begin{figure}[htbp]
    \centering
    \begin{minipage}{0.45\textwidth}
        \centering
        \caption{Number of data points and corresponding computation times for problem 2.}
        \begin{tabular}{|c|c|}
        \hline
        \textbf{\shortstack{Number of \\ Points}} & \textbf{\shortstack{ Computation Time \\ (seconds)}} \\
        \hline
        1000 & 14.855596 \\
        2000 & 25.246510 \\
        3000 & 37.736875 \\
        4000 & 49.941251 \\
        5000 & 64.811114 \\
        6000 & 74.934445 \\
        7000 & 87.540231 \\
        8000 & 98.864716 \\
        9000 & 112.871687 \\
        10000 & 124.236625 \\
        \hline
        \end{tabular}
        \label{table:computation_time_2}
    \end{minipage}
    \hfill
    \begin{minipage}{0.45\textwidth}
        \centering
        \caption{Plot of computation time vs. number of data points.}
        \includegraphics[width=\textwidth]{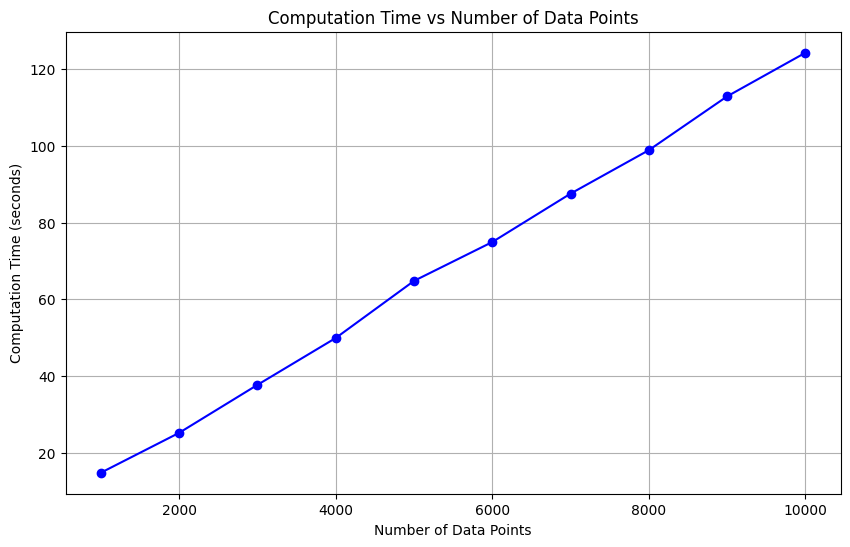}
        \label{fig:computation_time_plot_2}
    \end{minipage}
\end{figure}

\newpage
\section{Summary and conclusion}

In summary, this work enhances the PINN method by integrating the abstract mathematical structure of Lie groups. By altering the PINN model's loss function through infinitesimal transformations of the Lie group and introducing adaptivity in activation function, we achieve remarkable improvements. Numerical experiment on the quasi-linear Burgers equation, employing well reputed advanced activation functions, shows the effectiveness of our approach. Through the development of three distinct cases, we present the progressive evolution of PINN model, in case C we introduce m-ASPINN model. Our experiments prove that m-ASPINN is computationally efficient and attains comparable performance to state-of-the-art numerical methods. Consequently, m-ASPINN comes up with a promising solution, particularly for handling tasks requiring large number of datasets.
\\
In conclusion, the intersection of Lie symmetry groups and adaptivity into the PINN framework represents a significant betterment in computational modeling. The successful improvement and affirmation of the m-ASPINN model highlight its potential for addressing complex problems efficiently and accurately.

\bibliographystyle{unsrt} %
\bibliography{mybib}

\end{document}